\documentclass[12pt,psamsfonts]{amsart}

\newtheorem{anyprop}{Anyprop}[section]

\newtheorem{theorem}[anyprop]{Theorem}
\newtheorem{lemma}[anyprop]{Lemma}
\newtheorem{proposition}[anyprop]{Proposition}
\newtheorem{corollary}[anyprop]{Corollary}

\theoremstyle{definition}

\newtheorem{remark}[anyprop]{Remark}

\newcommand{\NN}{\mathbb{N}}

\newcommand{\PP}{\mathbb{P}}

\newcommand  {\shT}     {\mathcal{T}}

\newcommand  {\foa}     {\mathfrak{a}}

\newcommand  {\fom}     {\mathfrak{m}}

\newcommand  {\dual}    {\vee}

\newcommand  {\Ext}     {\operatorname{Ext}}

\newcommand  {\Grass}     {\operatorname{Grass}}

\newcommand  {\Hom}     {\operatorname{Hom}}

\newcommand  {\lra}     {\longrightarrow}

\renewcommand{\O}       {\mathcal{O}}

\newcommand  {\Pic}     {\operatorname{Pic}}

\newcommand  {\Proj}    {\operatorname{Proj}}

\newcommand  {\ra}      {\rightarrow}

\newcommand  {\rk}    {\operatorname{rk}}

\newcommand  {\Syz}     {\operatorname{Syz}}

\newcommand{\komdots}{ , \ldots , }

\newcommand{\oplusdots}{ \oplus \ldots \oplus }

\newcommand {\modu} {\, \rm mod \,}

\newcommand{\numiii}{\renewcommand{\labelenumi}{(\roman{enumi})}}

\newcommand{\stacklra}[1]{ \stackrel{ #1 }{\lra} }

\newcommand{\length} {\lambda}

\newcommand{\nuo}{\rho}
\newcommand{\nuu}{\sigma}
\newcommand{\nus}{\nu}

\newcommand{\pers }{\pi}

\newcommand{\hkf}{\varphi }

\newcommand{\degc}{\delta}
\newcommand{\expr}{\Phi}

\theoremstyle{remark}

\numberwithin{equation}{section}

\usepackage{amscd}
\usepackage{amssymb}

\def\mydate{\number\day\space\ifcase\month \or January\or February\or March\or April\or May\or
June\or July\or August\or September\or October\or November\or
December\fi \space\number\year}

\begin{document}

\title[Hilbert-Kunz function]
{Restriction of the cotangent bundle to elliptic curves and Hilbert-Kunz functions}

\author[Holger Brenner and Georg Hein]{Holger Brenner and Georg Hein}
\address{Mathematische Fakult\"at, Ruhr-Universit\"at Bochum, Universit\"atsstr 150,
               44780 Bochum, Germany;
Fachbereich Mathematik und Informatik,
Freie Universit\"at Berlin,
Arnimallee 3,
14195 Berlin, Germany.}

\email{Holger.Brenner@ruhr-uni-bochum.de, ghein@math.fu-berlin.de}

\begin{abstract}
We describe the possible restrictions of the cotangent bundle $\Omega_{\PP^N}$
to an elliptic curve $C \subset \PP^N$.
We apply this in positive characteristic
to the computation of the Hilbert-Kunz function
of a homogeneous $R_+$-primary ideal $I \subset R$ in the graded
section ring $R=\bigoplus_{n \in \NN} \Gamma(C, \O(n))$.
\end{abstract}

\maketitle


%===========================================================
\section*{Introduction}

The aim of this paper is to compute the Hilbert-Kunz function of homogeneous
primary ideals in the homogeneous coordinate ring
or in the graded  section ring of an elliptic curve $C \subset \PP^N$.
Recall that for a primary ideal $I$ in a local Noetherian ring $(R, \fom)$ or an $\NN$-graded algebra $(R,R_+)$
of dimension $d$ of positive characteristic $p$
the Hilbert-Kunz function is given
by
$$ e \mapsto \length ( R/I^{[p^{e}]}) \, ,$$
where $\length $ denotes the length.
The Hilbert-Kunz multiplicity of $I$ is defined as the limit
$$ e_{HK} (I)= \lim_{e \ra \infty}  \, \length (R/I^{[p^{e}]}) /p^{ed} \, .$$
This limit exists as a positive real number, as shown by Monsky in \cite{monskyhilbertkunz}.
For the maximal ideal $I= \fom$ this number is also
called the Hilbert-Kunz multiplicity of the ring $R$ itself.

It is an open question whether the Hilbert-Kunz multiplicity is always rational,
and the Hilbert-Kunz function is known only in very special cases.
The rationality in the case of a
homogeneous coordinate ring
$R=k[X_0 \komdots X_N]/\foa$ over a smooth projective curve
$C =V_+ (\foa) \subset \PP^N$ was
shown independently in \cite{brennerhilbertkunz} and \cite{trivedihilbertkunz}.
This work suggests a more geometric view on the problem of the Hilbert-Kunz function,
where stability properties of the syzygy bundle $\Syz(X_0 \komdots X_N)$ on the curve
are crucial. For example, if this bundle is strongly semistable,
meaning that every Frobenius pull-back of it is again semistable,
then the Hilbert-Kunz multiplicity of the cone is given by
$e_{HK}(R)=\deg (C) (N+1)/2N$ due to \cite[Corollary 2.7]{brennerhilbertkunz}.

The syzygy bundle $\Syz(X_0 \komdots X_N)$ on the curve is the restriction of the cotangent bundle $\Omega_{\PP^N} =\shT_{\PP^N}^\dual$, as follows from the Euler sequence.
What can we say about stability properties of the restriction of the (co)tangent bundle to a smooth curve?
This is still a difficult question in general, in particular in positive characteristic.
In this paper we investigate the case of cones over elliptic curves.
The main point which makes elliptic curves more accessible is the fact
that a semistable locally free sheaf on an elliptic curve is strongly semistable.

It is not difficult to show that the restriction of the cotangent bundle to an elliptic curve
embedded by a complete linear system is stable (Theorem \ref{tangentrestriction}), hence
indecomposable and strongly semistable.
We are however also interested in embedding linear systems which are not complete;
these linear systems exhibit a rich variety
of possible splitting behaviour of the restricted cotangent bundle.
We study in the first three sections
the vector bundles $E$ on an elliptic curve $C$
which appear as the restriction of the cotangent bundle $\Omega_{\PP^n}$ for some
embedding $\iota:C \to \PP^n$. This culminates
in the classification result Theorem \ref{spacecurves}, which holds in every characteristic
and is of independent interest.

In sections 4-6 we compute the Hilbert-Kunz function
of a homogeneous
$R_+$-primary ideal $I \subset R= \bigoplus_{n \geq 0} \Gamma(C, \O_C(n))$,
where $C \subset \PP^N$ is an elliptic curve over an algebraically closed field
of positive characteristic.
The general result is Theorem \ref{hilbertkunzfunctionelliptic},
which describes the Hilbert-Kunz function of $I=(f_1 \komdots f_n)$
as
$$ \hkf(q) = \frac{1}{2 \delta} ( \sum_{j=1}^l \frac{ \deg(S_j)^2}{\rk(S_j)} - \delta^2 \sum_{i=1}^n d_i) q^2
+ \gamma(q) \, ,$$
where $ \delta =\deg (\O_C)$, $\Syz(f_1 \komdots f_n)(0)= S_1 \oplusdots S_l$
is the decomposition into indecomposable sheaves,
$d_i=\deg(f_i)$ and $\gamma(q)$ is an eventually periodic function.

In section 5 and 6 we specialize this result to the computation of the Hilbert-Kunz function
of an ideal generated by the linear forms
$(X_0 \komdots X_N) \subseteq \Gamma(C, \O_C(1))$,
which define the embedding $C \subset \PP^N$.
If this ideal corresponds to a complete linear system
\--- that is $(X_0 \komdots X_N) = \Gamma(C, \O_C(1))$ \---,
then the semistability result for the restricted cotangent bundle mentioned above
yields the Hilbert-Kunz function of the normal
homogeneous coordinate rings over an elliptic curve.
This gives a new proof of a theorem of Fakhruddin-Trivedi (\cite{fakhruddintrivedi}),
which itself generalizes earlier results
of Monsky and Buchweitz-Chen in the case of a plane cubic curve (\cite{monskyhilbertkunz},
\cite{buchweitzchenhilbertkunz}).

In the final section 6 we deduce from the classification result for the restricted cotangent bundle
the possible Hilbert-Kunz functions for elliptic space curves $C \subset \PP^3$.

This paper was in part written when the first mentioned author was guest at the
Freie Universit\"at Berlin in february 2004. He would like to thank this institution for its hospitality. Furthermore we would like to thank the DFG-Schwerpunkt
``Globale Methoden in der komplexen Geometrie'' for financial support.

\section{The tangent bundle of projectively normal curves}

Throughout this section $C$, denotes a smooth projective
elliptic curve over an algebraically closed field $k$.
Recall that a vector bundle is called indecomposable if it is not the direct sum of two
proper subbundles.

The question we will discuss in the following sections is the following:
Consider an embedding of the elliptic curve $C$ into the projective space $\PP^r$.
This embedding is given by a surjection $\pi:\O_C ^{\oplus r+1} \to L$ where
$L$ is a very ample line bundle. The Euler sequence of $\PP^r$ gives that the kernel
of $\pi$ is the restriction of $\Omega_{\PP^r}(1)$ to the curve $C$.
We say that a rank $r$ bundle $F$ on $C$ is {\em a restricted twisted cotangent bundle},
if and only if there exists an embedding
$\iota : C \to \PP^r$ such that $F \cong \Omega_{\PP^r}(1)|_C$.
If $h^0(F) >0$, then there exists a linear subspace $\PP^{r-h^0(F)} \subset \PP^r$ such that $C$
is contained in this subspace.
Therefore, to classify all restricted twisted cotangent bundles
it is enough to assume $h^0(F)=0$.

In the following Proposition we list the facts about vector bundles on elliptic curves
we use thereafter. All these results can be found in Atiyah's article \cite{atiyahelliptic}.
In Tu's paper \cite{tusemistable}, these results are formulated in the frame work of semistable
bundles and their moduli spaces. See also \cite{heinploog} for a survey on these results
using Fourier--Mukai transforms.

\begin{proposition}
\label{ellcurbas}
(Properties of vector bundles on an elliptic curve $C$)
\vspace{-6pt}
\begin{itemize}
\item[(i)] Every vector bundle $E$ on $C$ decomposes into a direct sum of indecomposable ones
$E= \oplus E_i$.
\item[(ii)] There exist indecomposable vector bundles of degree $d$ and rank $r\geq 1$ on $C$.
These  vector bundles are stable, if and only if $d$ and $r$ are coprime.
\item[(iii)] An indecomposable vector bundle is semistable.
\item[(iv)] If $E$ is semistable of positive degree $d$, then we have $h^0(E)=d$, and $h^1(E)=0$.
This implies by Serre duality, that for $E$ semistable of negative degree, we have
$h^0(E)=0$ and $h^1(E)=- \deg(E)$.
\item[(v)] If $E$ is semistable of degree zero, then we consider the following set:
$\Theta_E:= \{ L \in \Pic^0(C) \, | \, h^0(E \otimes L) \ne 0\}$. The subset
$\Theta_E$ of $\Pic^0(C)$ is non-empty with at most $\rk(E)$ elements.
If $E$ is furthermore indecomposable, then $\Theta_E$ has exactly one element.
\end{itemize}
\end{proposition}

\begin{lemma}
\label{globgen}
Let $E$ be a semistable vector bundle on the elliptic curve $C$.
If $E$ is not isomorphic to the bundle $\O_C^{\oplus r}$, then the following holds:
$E$ is globally generated, if and only if $\mu(E) >1$.
\end{lemma}
\proof
In order to study global generatedness, we consider the following long exact sequence
for all points $P \in C$:
$$H^0(E) \stackrel{\psi_P}{\longrightarrow} E \otimes k(P) \to H^1(E(-P)) \to H^1(E) \to 0\,.$$
If $\mu(E)>1$, then $\mu(E(-P))>0$. Thus, by proposition \ref{ellcurbas} (iv),
we have $H^1(E(-P))=0$
for all points $P$. Thus, all maps $\psi_P$ are surjective.

Suppose now, that $E$ is globally generated. In other words, $E$ is a quotient of the semistable
bundle $H^0(E) \otimes \O_C$. Thus, $\mu(E) \geq 0$ with equality only for
$E \cong H^0(E) \otimes \O_C$. Since this case was excluded, we have $H^1(E)=0$.
Since all the maps $\psi_P$ are surjective, we conclude from the above exact
sequence that $H^1(E(-P))=0$ for all points $P \in C$.
Since all line bundles in $\Pic^0(C)$ are of type $\O_C(P_0-P)$, we conclude from part (v) of
proposition \ref{ellcurbas} that $\mu(E(-P_0)) > 0$.
Consequently, $\mu(E) >1$.
\qed

This lemma yields the stability of the restricted tangent bundle for any
embedding of $C$ be a complete linear system:

\begin{theorem}
\label{tangentrestriction}
Let $L$ be a line bundle on $C$ of degree $d \geq 2$.
The kernel of the evaluation map $H^0(L) \otimes \O_C \to L$ is stable of
degree $-d$ and rank $d-1$.
\end{theorem}
\proof
Denote the kernel of the evaluation map with $K$. By construction, we have $h^0(K)=0$.
Since the rank and the degree of $K$ are coprime, semistability implies stability.
Assume that $K$ is not semistable. Let $K_1$ be a subsheaf of maximal slope. This implies that
$K_1$ is semistable. On the one hand, we have $h^0(K_1)=0$,
because $K_1$ is a subbundle of $K$. 
On the other hand, we have $\mu(K_1) > \mu(K) = -1 -\frac{1}{d-1}$.
Since the rank of $K_1$ is less then $d-1$, we conclude $\mu(K_1) \geq -1$.

Dualizing the embedding $K_1 \to K \to H^0(L) \otimes \O_C$,
we obtain that $E:=K_1^\lor$ is globally generated, and by Serre duality $h^1(E)=0$.
However, since $\mu(E) = -\mu(K_1) \leq 1$, this contradicts the assertion of lemma
\ref{globgen}.
\qed

\section{Embedding vector bundles}

Let $E$ be a vector bundle on $C$ which is globally generated and satisfies $h^1(E)=0$,
then $V:=H^0(E)$ is a $k$-vector space of dimension $d = \deg(E)$.
The surjection $p:V \otimes \O_C \to E$ defines a morphism $\tilde p: C \to \Grass(V,r)$,
where $\Grass(V,r)$ is the Grassmannian variety of $r$ dimensional quotients of $V$.
We next investigate when this morphism is an embedding. To do so, we consider the decompostion
of $E$ into indecomposable bundles $E= \oplus E_i$.
We consider now the decompostion of $E=E_{\mu <2} \oplus E_{\mu=2} \oplus E_{\mu > 2}$ 
with the three summands defined by
$$ E_{\mu <2}:= \bigoplus_{\mu(E_i) <2}  E_i\, , \qquad
E_{\mu =2}:= \bigoplus_{\mu(E_i) =2}  E_i\, , \qquad \mbox{and }
E_{\mu >2}:= \bigoplus_{\mu(E_i) >2}  E_i\,.  $$
The next lemma investigates, whether $\tilde p$ is
an embedding.

\begin{lemma}
Let $E$ be a globally generated vector bundle of degree $d$ and rank $r$ with $h^1(E)=0$.
We denote by $V$ its space of global sections.
The morphism $\tilde p:C \to \Grass(V,r)$ corresponding to the surjection
$p:V \otimes \O_C \to E$ is an embedding unless $E = E_{\mu=2}$, and $E_{\mu=2}$ is
isomorphic to $L^{\oplus \rk(E)}$ for a line bundle $L \in \Pic^2(C)$.
\end{lemma}

\proof
Since $E$ is globally generated, we have $h^0(E(-P))=d-r$.
Analogously to \cite[Proposition IV.3.1]{haralg}, we see that $\tilde p$ is an
embedding, if and only if for every two points $P,Q \in C$
(including the case $P=Q$) the inequality
$h^0(E(-P-Q)) < d-r$ holds.

Suppose that $\tilde p$ does not define an embedding. We conclude that there exists two point
$P,Q \in C$ such that $h^0(E(-P-Q))=h^0(E(-P))$. This implies that we have three equalities
$h^0(E_{\mu<2}(-P-Q))=h^0(E_{\mu<2}(-P))$, $h^0(E_{\mu=2}(-P-Q))=h^0(E_{\mu=2}(-P))$, and
$h^0(E_{\mu>2}(-P-Q))=h^0(E_{\mu>2}(-P))$. However by proposition \ref{ellcurbas} (iv),
we conclude
that $h^0(E_{\mu>2}(-P-Q))=h^0(E_{\mu>2}(-P)) - \rk(E_{\mu > 2})$. This yields $\rk(E_{\mu > 2})=0$.
Hence, $E_{\mu > 2}=0$. Analogously we deduce that $E_{\mu < 2}=0$.

This yields that $E = E_{\mu=2}$. Consquently we have $\deg(E)=2r$, and 
$h^0(E(-P)=r$. If $h^0(E(-P-Q))=r$, then we have a semistable sheaf of degree zero
with $r$ global sections. This implies $E(-P-Q) \cong \O_C^{\oplus r}$.
Thus, $E \cong \left( \O_C(P+Q) \right)^{\oplus r}$, as stated.

On the other hand, we see that for $E \cong \left( \O_C(P+Q) \right)^{\oplus r}$
the corresponding morphism $\tilde p$ is $2 : 1$ onto its image. Reversing the above
arguments, we see that $\tilde p$ defines an embedding otherwise.
\qed

Suppose now that $\tilde p: C \to \Grass(V,r)$ defines an embedding.
If this is the case, then we will call $E$ an {\em embedding vector bundle\/}.
If the dimension of $V$ equals $r+1$, then we are in the dual situation of 
theorem \ref{tangentrestriction}. If $W \subset V$ is a subspace of
dimension $r+1$, then we obtain a rational map $\tilde p_W:C - - > \Grass(W,r)$.

We decide next whether $\tilde p_W$ is a regular morphism to, or even
an embedding into $\Grass(W,r)$.
We need the following dimension lemma:

\begin{lemma}
\label{dimcycle}
Let $V$ be a $k$-vector space of dimension $d$.
Moreover, fix integers $r$, $l$, and $k$ satisfying the inequalities
$$1 < r < d-2 \,, \quad 0 \leq k \leq d-r \,, \quad \mbox{and} \quad l \geq 1\,.$$
For a subspace $V'$ of $V$ of codimension $r+k$,
we consider the following subset $Z$ of the Grassmannian $\Grass(r+1,V)$ of $r+1$
dimensional subspaces of $V$:
$$Z:= \{ W \subset V | \, \dim(W)=r+1 \,, \quad
\mbox{and} \quad \dim(W \cap V') \geq l \} \,.$$ 
This $Z$ is a closed subset of $\Grass(r+1,V)$ of codimension $l^2+lk-l$.
\end{lemma}
\proof This is a standard dimension computation and can be done as follows.
Consider the projection $\pi:V \to V'':=V/V'$. First, we choose a
$(r+1-l)$-dimensional subspace $W''$ of $V''$. Secondly, we choose $W$
to be an $(r+1)$-dimensional subspace of $\pi^{-1}(W'')$.
Using the fact that the dimension of the Grassmannian of $m$ dimensional subspaces
of an $n$ dimensional space is $m(n-m)$, the lemma follows immediately.
\qed

As a consequence, we obtain the following lemma, telling us that most subspaces
$W$ of $V$ define a regular morphism $\tilde p_W$.

\begin{corollary}
\label{general-gen}
If a rank $r$ vector bundle $E$
on a curve $C$ is globally generated with $h^0(E) > r+1$,
then the $r+1$ dimensional subspaces $W$ of $H^0(E)$,
which don't generate $E$,
form a proper closed subset of the Grassmannian $\Grass(r+1,H^0(E))$.
\end{corollary}
\proof For a point $P \in C$ we have that $H^0(E(-P))$ is of codimension
$r$ in $H^0(E)$. If $W$ does not generate $E$ at $P$, then we must have
$\dim(W \cap H^0(E(-))) \geq 2$. By lemma \ref{dimcycle} these $W$ form
a closed subscheme of codimension $2$. Considering all points $P$ of $C$,
we obtain a
one dimensional family of codimension two subschemes of $\Grass(r+1,H^0(E))$.
\qed

For a globally generated sheaf $E$ on the elliptic curve $C$, we define the following two numbers:
$$a_E := \min_{P,Q \in C} h^0(E(-P-Q)) \,, \quad \mbox{and} \quad
b_E:= \max_{P,Q \in C} h^0(E(-P-Q)) \,.$$

\begin{lemma}
For a globally generated vector bundle $E$ of rank $r$ on an elliptic curve $C$ we 
have:
\vspace{-6pt}
\begin{itemize}
\item[(i)] $h^0(E)-2 \cdot r \leq a_E \leq b_E \leq h^0(E)-r$;
\item[(ii)] The set $U$ of all pairs $(P,Q) \in C \times C$ with $h^0(E(-P-Q))=a_E$ is Zariski open
and nonempty.
\item[(iii)] If $a_E < b_E$,
then the complement of $U$ in $C \times C$ consists of a finite union of fibers of the
addition map $+ : C \times C \to C$.
\end{itemize}
\end{lemma}
\proof
(i) is trivial.
(ii) is just the semicontinuity theorem (see \cite[Theorem III.12.8]{haralg}).
To prove (iii) we note that for $P+Q = P'+Q'$, we have $\O_C(-P-Q) \cong \O_C(-P'-Q')$.
Thus, $E(-P-Q) \cong E(-P'-Q')$. However note that $E(-P-Q)$ and $E(-P'-Q')$ are
different subsheaves of $E$.
\qed

We use these numbers to define two conditions:

We say that a globally generated vector bundle $E$ of rank $r$ satisfies\\
condition $(*)$,
if $a_E=0$, or $a_E < h^0(E)-r-2$ holds, and\\
condition $(**)$, if $b_E=0$, or $b_E < h^0(E)-r-1$ holds.

Using these two conditions and the dimension formula from lemma \ref{dimcycle} we obtain the
following:
\begin{corollary}\label{sufficient}
Let $E$ be an embedding vector bundle of rank $r \geq 3$ with $h^1(E)=0$ on an elliptic curve $C$.
If $E$ satisfies the conditions $(*)$ and $(**)$, then the general $(r+1)$-dimensional subspace
$W$ of $H^0(E)$ is associated to an embedding $\tilde p_W: C \to \PP(W^\lor)$.
\end{corollary}
\proof
We have seen in corollary \ref{general-gen} that a general subspace $W$ of dimension $r+1$ in
$H^0(E)$ generates $E$. If $W$ does not define an embedding, there must be two points $P$ and
$Q$ in $C$ such that $W \cap H^0(E(-P)) = W \cap H^0(E(-Q))$. This implies that $W \cap H^0(E(-P-Q))$
is of dimension one. Now we see that lemma \ref{dimcycle} and the conditions $(*)$, and $(**)$
yield that those $W$, for which there exists a pair $(P,Q) \in C^2$ with $W \cap H^0(E(-P-Q))$
is not zero, form a proper closed subscheme of $\Grass(r+1,H^0(E))$.
\qed

\begin{lemma}
\label{nostar}
Let $E$ be an embedding vector bundle of rank $\rk(E) \geq 3$ with $h^1(E)=0$.
If $E$ does not satisfy the condition $(*)$, then $E$ is isomorphic to the direct sum
$E_1 \oplus L$, where $E_1$ is a stable vector bundle of rank $\rk(E)-1$ and degree
$\rk(E)$, and $L$ is a line bundle of degree $\deg(E)-\rk(E) \geq 3$. 
\end{lemma}
\proof
We decompose $E$ into the direct sum $E=E_{\mu<2} \oplus E_{\mu=2} \oplus E_{\mu>2}$.
Since by parts (iv) and (v) of proposition \ref{ellcurbas} we have
$h^0(E(-P-Q)) = h^0(E_{\mu>2}(-P-Q))$ for
two general points $P$ and $Q$, we deduce from $a_E >0$ that $E_{\mu > 2} \ne 0$.
Since $E$ is an embedding bundle, we must have $a_E< h^0(E)-\rk(E)$. To violate condition $(*)$
we have to consider the following two cases:

{\em Case 1:} $a_E = h^0(E)-\rk(E)-1$.\\
If $E_{\mu<2}$ were trivial, then we would have $a_E=h^0(E)-2\rk(E)$ which would imply $E$ is
a line bundle. Thus, $E_{\mu<2}$ is nontrivial. From lemma \ref{globgen} we conclude
$\deg(E_{\mu<2})>\rk(E_{\mu<2})$. We have now
$1=(\deg(E_{\mu<2})-\rk(E_{\mu<2}))+ \rk(E_{mu=2}) + \rk(E_{mu >2})$.
In this sum of nonnegative integers the first and the last summand are
positive which is a contradiction.
Thus, the only possibility is the remaining

{\em Case 2:} $a_E = h^0(E)-\rk(E)-2$.\\
As before, we deduce the equality $2=(\deg(E_{\mu<2})-\rk(E_{\mu<2}))+ \rk(E_{mu=2}) + \rk(E_{mu >2})$,
with the first and last summand positive. Hence, we must have $\deg(E_{\mu<2})-\rk(E_{\mu<2})=1$, and
$E_{\mu >2}$ is a line bundle.

If $E_{\mu<2}=\bigoplus_{i=1}^lE_i$ is the decomposition into its indecomposable components,
the $\deg(E_i)> \rk(E_i)$ has to hold for each component which implies immediately
that $E_{\mu<2}$ is indecomposable of rank $\rk(E)-1$ and degree $\rk(E)$.
\qed

\begin{lemma}
\label{nostars}
Let $E$ be an embedding vector bundle of rank $r=\rk(E) \geq 3$ with $h^1(E)=0$.
If $E$ does not satisfy the condition $(**)$, then $E$ is isomorphic to
one of the following vector bundles

\vspace{-6pt}
\begin{itemize}
\item[(i)] $E_1 \oplus L^{\oplus r_2}$\\
$E_1$ is stable of rank $(\rk(E)-r_2)>2$,
$r_2>0$, $\deg(E_1)=\rk(E_1)+1$, and $L \in \Pic^2(C)$;
\item[(ii)] $L^{\oplus(r-1)} \oplus M$\\
$L$ and $M$ are line bundles with
$\deg(L)=2$ and $\deg(M) \geq 3$;
\item[(iii)] $L^{\oplus(r-1)} \oplus M$\\
$L$ and $M$ are line bundles of degree two with $L \ncong M$;
\item[(iv)] $L^{\oplus(r-2)} \oplus (L \otimes E_2)$\\
$L$ is a line bundle of degree two, and $E$ is the
nontrivial extension bundle in $\Ext^1(\O_C,\O_C)$.
\end{itemize}
\end{lemma}
\proof
Since $E$ is an embedding bundle, we have $b_E < h^0(E)-r$.
Thus, to violate the condition $(**)$, we must have $0 < b_E =h^0(E)-r-1$.
We consider the direct sum composition $E=E_{\mu <2} \oplus E_{\mu=2} \oplus E_{\mu > 2}$.
Analogously to the proof of lemma \ref{nostar}, we conclude that
$$b_E \leq h^0(E)-r-(\deg(E_{\mu <2})-\rk(E_{\mu < 2})) - \rk(E_{\mu >2})\, .$$
Hence, at least one of the bundles $E_{\mu <2}$ or $E_{\mu>2}$ is trivial.
Let $P$ and $Q$ be two points of the curve $C$ such that $h^0(E(-P-Q))=b_E$.
We distinguish the three cases:\\
{\em Case 1: $E_{\mu <2} \ncong 0$.}\\
As in case 2 of lemma \ref{nostar}, we conclude that $E_{\mu < 2}$ is stable of rank $r_1$ and of
degree $r_1+1$. Thus, $E_{\mu =2}$ turns out to be a bundle of rank $r_2=r-r_1$, which
satisfies $h^0(E_{\mu = 2}(-P -Q))=r_2$. Knowing that $E_{\mu =2}$ is semistable of
slope zero with $r_2$ global sections, we deduce that $E_{\mu=2}(-P-Q) \cong \O_C^{\oplus r_2}$.
Thus, we are in case (i) of the lemma.\\
{\em Case 2: $E_{\mu >2} \ncong 0$.}\\
From the above inequality, we deduce that $E_{\mu >2}$ is a line
bundle and, as in case 1, we obtain $h^0(E_{\mu =2}(-P-Q)) = \rk(E_{\mu =2})$.
Thus, the assumption of case 2 yields case (ii) of the lemma.\\
{\em Case 3: $E = E_{\mu =2}$.}\\
From $h^0(E(-P-Q))=r-1$, we deduce that $E(-P-Q)$ contains $\O_C^{\oplus r-1}$ as a subbundle.
We consider the short exact sequence
$$0 \to \O_C^{\oplus r-1} \to E(-P-Q) \to M \to 0 \, .$$
$M$ is a line bundle of degree zero. If $M$ is not isomorphic to $\O_C$, then $\Ext^1(M,\O_C)=0$
and this sequence splits. Hence $M \ncong \O_C$ implies that we are in case (iii).
If $M \cong \O_C$, then the above short exact sequence has to be non split because
$h^0(E)=r-1$. Thus, we are in the remaining case (iv) of our lemma.
\qed

\begin{lemma}\label{necessary1}
Let $E$ be an embedding vector bundle with $h^1(E)=0$ of rank $r \geq 3$.
If $E$ does not satisfy condition $(*)$, then for all surjections $p: \O_C \otimes W \to E$,
where $W$ is a vector space of dimension $r+1$,
the corresponding morphism $\tilde p:C \to \Grass(W,r)$ is not an embedding.
\end{lemma}
\proof By lemma \ref{nostar},
we know that $E$ is isomorphic to $E_1 \oplus L$ with $L$ a line bundle,
and $E_1$ stable of degree $r$. Suppose that $p: \O_C \otimes W \to E$ is
surjective for $W$ of
dimension $r+1$.
Consider the composition
$\O_C \otimes W \to E \to E_1$. Since $h^0(E_1)=r$, we deduce that 
there exists a nontrivial section $s$ of $L$ in $W$.
Since $\deg(L) \geq 3$, we have two points $P$ and $Q$ in $C$ 
(possibly $P=Q$) such that $s \in H^0(L(-P-Q))$. However, this means $W \cap H^0(E(-P-Q))$
is of dimension one. Thus, $\tilde p$ does not separate $P$ and $Q$.
\qed

\begin{lemma}\label{necessary2}
Let $E$ be an embedding vector bundle with $h^1(E)=0$ of rank $r \geq 3$.
If $E$ does not satisfy condition $(**)$,
then for all surjections $p: \O_C \otimes W \to E$,
where $W$ is a vector space of dimension $r+1$,
the corresponding morphism $\tilde p:C \to \Grass(W,r)$ is not an embedding.
\end{lemma}
\proof
We need to examine the four cases of lemma \ref{nostars}.
We keep the notation of that lemma, and consider an arbitrary surjection
$p:\O_C \otimes W \to E$ where $W$ is a vector space of dimension $r+1$.
We consider $W$ as a subspace of $H^0(E)$.
In order to show that $\tilde p$ is not an embedding,
we have to find two points $P$ and $Q$ in $C$
such that $H^0(E(-P-Q)) \cap W$ is one-dimensional. This implies that $\tilde p$ does not
separate these two points.

{\em Case (i): $E \cong E_1 \oplus L^{\oplus r_2}$.}\\
Consider the resulting morphism $p_1: W \to H^0(E_1)$.
Its kernel is a subspace of dimension at least $r_2$.
Take $W' \subset W$ to be an $r_2$ dimensional subspace of $W$ in the kernel
of $p_1$.
We have $W' \subset H^0(L^{\oplus r_2})$. The corresponding morphism
$q:\O_C \otimes W' \to L^{\oplus r_2}$ cannot be surjective.
Thus, there exists a point $P \in C$ and a nontrivial section $s \in W'$
with $s \in H^0(L^{\oplus r_2}(-P))$.
There exists a unique point $Q \in C$ such that $L(-P) \cong \O_C(Q)$.
This implies that $s$ also lies in $H^0(L^{\oplus r_2}(-P-Q))$.
Thus, $s \in W \cap H^0(E(-P-Q))$.

{\em Case (ii): $E \cong L^{\oplus r-1} \oplus M$.}\\
Let $L \cong \O_C(Q_1+Q_2)$. For a point $P \in C$, there exists a unique point $\psi(P)$
such that $L \cong \O_C(P + \psi(P))$. Indeed, putting $\psi(P)=Q_1+Q_2-P$,
we obtain this point.
For a point $P$, we consider the kernel of the map $W \to E \otimes k(P)$.
Since $E$ is globally generated by $W$, this is spanned by a global section $s$ of $E(-P)$. 
The global section $s$ yields two global sections $s_1 \in H^0(L(-P)^{\oplus r-1})$ and section
$s_2 \in H^0(M(-P))$. Since $h^0(L(-P))=1$, the section $s_1$ vanishes in $\psi(P)$, too.
Thus, to show that $\O_C \otimes W \to E$ does not define an embedding, it suffices to show
that there exists at least one point $P$ such that the global section $s_2 \in H^0(M(-P))$
also vanishes at $\psi(P)$.
To do so, we consider the two projections $p$ and $q$ from
$C \times C$ to the components
$$C \stackrel{q}{\leftarrow} C \times C \stackrel{p}{\rightarrow} C$$
From the short exact sequence
$0 \to K \to \O_C \otimes W \to E \to 0$ on $C$,
we construct the following diagram
on $C \times C$
$$\begin{array}{ccccccc}
0 \to & p^*K & \to & W' & \to & p^*E(-\Delta) & \to 0 \\
& || && \downarrow && \downarrow\\
0 \to & p^*K & \to & \O_{C \times C} \otimes W & \to & p^*E & \to 0\\
&&& \downarrow && \downarrow\\
&&& p^*E|_\Delta & = &  p^*E|_\Delta\\
\end{array}$$
with $\Delta$ denoting the diagonal in the product $C \times C$.
Applying the functor $q_*$ to the middle column,
we obtain that $q_*W' \cong K$.
We obtain a morphism $q^*K \to q_*(p^*M(-\Delta))$. This morphism assigns in any point $P$ of $C$
the global section $s_2$ of $M(-P)$.
To decide whether this section vanishes at $\psi(P)$,
we just have to intersect the divisors
$\psi(\Delta)=\{ (P,\psi(P)) \, | \, P \in C \}$ with the vanishing divisor $D$ of
$q^*K \to p^*M(-\Delta)$. Up to numerical equivalence, we obtain
$D=\deg(M)F_p + (\deg(M)+2r-2)F_q -\Delta$, where $F_p$ (resp. $F_q$) denote a fiber of $p$ (resp. $q$).
Since $D.\psi(\Delta)= 2 \deg(M) + 2 r -6 >0$,
we have  points $P$ such that $s_2$ vanishes in $\psi(P)$ as well, which
contradicts the embedding.

{\em Case (iii) and (iv): There exist two line bundles $L$ and $M$ of degree two,
and a short exact sequence
$0 \to L^{\oplus r-1} \to E \to M \to 0$.}\\
Considering the kernel of $W \to H^0(M)$ we find again an $r-1$ dimensional
subspace $W'$ of
global sections of $L^{\oplus r-1}$ in $W$. Now we proceed as in case (i).
\qed

\section{Classifying the restriction of the cotangent bundles to elliptic curves}

Now we have everything at hand to give a classification of all restricted twisted cotangent
bundles. We remark that such a bundle $F$ is the direct sum $\O_C^{\oplus u} \oplus F'$
where $F'$ is a restricted twisted cotangent
bundle  without global sections.
Therefore, we content ourselves with the classification of all
those without global sections.

\begin{theorem}
\label{spacecurves}
Let $F$ be a vector bundle on an elliptic curve $C$
of rank $r \geq 3$. Assume furthermore that $h^0(F)=0$. 
Let $F=\bigoplus_{i=1}^lF_i$ be the decomposition of $F$ into its indecomposable 
components. Assume that we have $\mu(F_1) \geq \mu(F_2) \geq \ldots \geq \mu(F_l)$.
The bundle $F$ is a restricted twisted cotangent bundle, if and only if the
conditions (i)--(iv) hold.

\vspace{-6pt}
\begin{itemize}
\item[(i)] The biggest slope $\mu(F_1)$ satisfies $\mu(F_1) < -1$;
\item[(ii)] $F$ is not isomorphic to the direct sum $F_1 \oplus M$ with
$M$ a line bundle and $F_1$ stable of degree $-r$.
\item[(iii)] $F$ is not isomorphic to $F_1 \oplus L^{\oplus r_2}$ with
$0 < r_2 \leq r-2$, $\deg(L)=-2$, and $F_1$ stable of degree $r_2-r-1$. 
\item[(iv)] There exists no exact sequence $0 \to L^{\oplus r-1} \to F \to M \to 0$
where $L$ is a line bundle of degree $-2$ and $M$ a line bundle of degree $\deg(M) \leq -2$.
\end{itemize}
\end{theorem}
\proof
If $\iota:C \to \PP^n$ is an embedding, then the pull back of the Euler sequence
yields
$0 \to \iota^*\Omega_{\PP^n}(1) \to \O_C \otimes V \to \O_C(1) \to 0$. The morphism
$\iota$ is dual to the morphism $\iota^\lor:C \to \Grass(V^\lor,n)$
defined by the surjection $\O_C \otimes V^\lor \to (\iota^*\Omega_{\PP^n}(1))^\lor$.
Thus, $\iota^\lor$ is an embedding too. Consequently we have that the bundle
$E=(\iota^*\Omega_{\PP^n}(1))^\lor$ is an embedding vector bundle which 
defines via the surjection of a vector space of dimension $\rk(E)+1$ an embedding
into projective space. Now the result follows from the sufficient (Corollary
\ref{sufficient}) and necessary (Lemma \ref{necessary1} and Lemma \ref{necessary2})
conditions $(*)$ and $(**)$ and the classification of all globally generated
vector bundles which do not satisfy these conditions (Lemma \ref{nostar} and Lemma
\ref{nostars}).
\qed

Let us classify all rank three vector bundles which are restricted twisted
cotangent bundles. There are only three cases for the ranks of the decomposition
$F=\bigoplus_{i=1}^lF_i$ into indecomposable bundles.

{\em Case (i): $l=1$, i.e. $F$ is indecomposable.}\\
$F$ is a restricted twisted cotangent bundle, if and only if $\deg(F) \leq -4$.

{\em Case (ii): $l=2$, i.e. $F$ is the direct sum of a line bundle $L$ and a rank
two vector bundle $F_1$.}\\
$F$ is a restricted cotangent bundle if $\deg(F_1) \leq -4$, and $\deg(L) \leq -2$.
In the case that in both inequalities the equality holds, we must furthermore demand that
$\Hom(F_1,L)=0$.

{\em Case (iii): $l=3$, i.e. $F$ is a direct sum of three line bundles $L_1 \oplus L_2 \oplus L_3$.}\\
$F$ is a restricted cotangent bundle if $\deg(L_i) \leq -2$,
and any two of these line bundles which are
of degree $-2$, are not isomorphic to each other.

\section{The Hilbert-Kunz function over an elliptic curve}
\label{functiongeneralsection}

In this section we compute the Hilbert-Kunz function of an $R_+$-primary homogeneous ideal
$I \subseteq R$,
where $R= \bigoplus _{n \in \NN} \Gamma(C,\O(n))$ is the graded section ring associated
to a very ample invertible sheaf $\O(1)$ on an elliptic curve $C$.
Homogeneous ideal generators $I=(f_1 \komdots f_n)$ of degree $d_i = \deg (f_i)$
give rise to the short exact sequence on $C$,
$$0 \lra \Syz(f_1 \komdots f_n)(m) \lra \bigoplus_{i=1}^n \O(m - d_i) \stackrel{f_i}{\lra} \O(m) \lra 0 \, ,$$
which relates the computation of the Hilbert-Kunz function
to the computation of the global sections of
$\Syz(f_1 \komdots f_n)(m)$
and of its Frobenius pull-backs. We first need the following Lemma.

\begin{lemma}
\label{indecomposablesections}
Let $C$ denote an elliptic curve over an algebraically closed field $k$
of positive characteristic $p$
together with a fixed very ample invertible sheaf $\O(1)$ of degree $\degc = \deg(C)$.
Let $S$ denote an indecomposable sheaf on $C$.
Set $\nus = - \mu(S) /\degc = - \deg(S)/ \rk(S) \delta $.
Let $\nuu$ and $\nuo$ denote natural numbers such that
$ \nuu \leq \nus < \nuo$.
Write
$ \lceil q \nus \rceil =  q \nus + \pers (q) $
with an eventually periodic function $\pers (q)$. Then
\begin{eqnarray*}
\sum_{m=q \nuu }^{q \nuo  -1}  h^0(S^q(m)) 
& =&  \frac{q^2}{2 \degc} 
\big( \frac{\deg(S)^2 }{\rk(S)} +  2 \nuo \deg(S) \degc + \nuo^2 \rk(S) \degc ^2 \big) \cr
& & - q \frac{\degc}{2}  (\nuo \rk(S) + \frac{\deg(S)}{\degc})  
 +   \pers (q)(1 - \pers (q) )\rk(S)   \frac{\degc }{2}  \cr
& & + h^1(S^q(   \lceil q \nus \rceil))
\end{eqnarray*}
\end{lemma}
\proof
An indecomposabe sheaf on an elliptic curve is strongly semistable.
For $m < \lceil q \nus \rceil $ we have
$m< -q \mu(S) / \degc$ and therefore
$$ \deg(S^q(m)) = q \deg(S) + m \rk(S)  \degc <  0\, .$$
Since $S^q$ is semistable we have $h^0(S^q(m))=0$ for $ m < \lceil q \nus \rceil $.
So the sum starts with $m = \lceil q \nus \rceil $.

For $m > \lceil q \nus \rceil $ we find that
$h^1( S^q(m))=0$ by semistability and Serre duality.
Hence by Riemann-Roch we get
\begin{eqnarray*}
\sum_{m = \lceil q \nus \rceil}^{q \nuo -1 } h^0(S^q(m)) &  = & 
\sum_{m = \lceil q \nus \rceil}^{q \nuo -1} \deg(S^q(m))
+ h^1(S^q(\lceil q \nus \rceil)) \cr
&=& \sum_{m = \lceil q \nus \rceil}^{q \nuo -1 }  \big( q\deg(S) +m\rk(S) \degc \big)
 + h^1(S^q(\lceil q \nus \rceil))      \cr
&=& q ({ q \nuo } -{\lceil q \nus \rceil}) \deg(S)  \cr
& &
+\big(  q \nuo  (  q \nuo -1) - \lceil q \nus \rceil (\lceil q \nus \rceil - 1) \big)
\frac{ \rk(S)\degc}{2} \cr
& &   +h^1(S^q(\lceil q \nus \rceil))    \, .
\end{eqnarray*}
We write $\lceil q \nus \rceil =  q \nus + \pers (q)$. The first summand is
$(q^2 \nuo - q^2 \nus -q \pers (q) ) \deg(S)$
and the second summand is
$$ \big( q^2 (\nuo^2-\nus^2) -q(\nuo - \nus) - 2 q \nus \pers (q) -\pers (q)(\pers (q)-1) \big)
\frac{ \rk(S)\degc}{2} \, .$$
We regroup and use $\nus =- \deg (S)/ \rk (S) \degc$
to get the stated result.
\qed

\begin{remark}
\label{periodicityremark}
Let $S$ denote an indecomposable bundle as in lemma \ref{indecomposablesections},
set $\nu= - \deg(S)/ \rk(S) \degc$. Then we have defined $\pers (q)$
by the condition $\lceil q \nu \rceil =q \nu + \pers (q)$.
Since $\lceil a/b \rceil = a/b + ( -a \modu b )/b$,
we have
$$\pers (q) =  ( q\deg (S) \modu \rk(S) \degc) / \rk(S) \degc \, .$$
This is eventually a periodic function.
In important cases, e.g. if $S$ is an indecomposable syzygy bundle,
the degree $\deg(S)$ is a multiple of the curve degree, say
$\deg(S)=u \degc$; then we can write $\pers(q)=( qu \modu \rk(S))/ \rk(S)$
and the length of the periodicity is bounded by $\rk(S)$.
\end{remark}

\begin{remark}
\label{h1remark}
Suppose again the situation of the previous lemma.
If $q \nu = -q \deg (S) / \rk(S) \degc$ is not an integer,
then $\deg(S^q(\lceil q \nu \rceil))>0$ and the $h^1$-term vanishes.
So suppose that $q \nu$ is an integer so that
$S^q(q \nu)$ has degree $0$.
Due to
$S^{qp} (qp\nu) =(S^q(q \nu))^p$ we see that
$\deg ( S^q(q \nu))=0$ for all $q \gg 0$.
Since $S^q(q \nu)$ is strongly semistable, we have a decomposition
$S^q(q \nu)= G_1 \oplusdots G_t$
with indecomposable bundles of degree $0$.
Every $G_j$ equals $G_j =F \otimes L$, where $L$ is an invertible sheaf of degree $0$
and where $F$ is the unique indecomposable sheaf of degree $0$ and rank $\rk(F)= \rk (G_j)$
with a non-trivial global section (see proposition \ref{ellcurbas} (v)).

If the Hasse-invariant of the elliptic curve is $0$, then $F^q$ is trivial
for $q \gg 0$. In this case $S^q( q \nu)$ is a direct sum of line bundles of degree $0$, $q \gg 0$.
For a line bundle $L$ either $L^q = \O$ for some $q$ or not.
Therefore for $q \gg 0$ we have $S^q(q \nu)= \O^u \oplus L_1 \oplusdots L_m$, where
$L_j$ are line bundles such that no $q$-power of them gets trivial.
The number of sections of $S^q(q \nu)$ as well as the dimension
of $H^1(S^q(\lceil q \nu \rceil))$
is then eventually constant $=u$.

If the Hasse-invariant of the elliptic curve is $1$,
then there exists also a decomposition
$S^q(q \nu)= \O^u \oplus G_1 \oplusdots G_t$ such that the $G_j$ are strongly stable.
If $G_j =F \otimes L$, then $(G_j)^p =F \otimes L^p$, and the dimension of
the global sections of $S^q(q \nu)$ is again eventually constant.
\end{remark}

\begin{theorem}
\label{hilbertkunzfunctionelliptic}
Let $k$ denote an algebraically closed field of positive characteristic $p$.
Let $C \subset \PP^N$ denote an elliptic curve
of degree $\degc = \deg (C) = \deg (\O(1)|C)$.
Set
$R= \bigoplus_{n \in \NN} \Gamma(C, \O(n))$.
Let $I$ denote a homogeneous
$R_+$-primary ideal in $R$ given by homogeneous
ideal generators $I=(f_1 \komdots f_n)$ of degree $d_i =\deg(f_i)$.
Let $\Syz(f_1 \komdots f_n)(0) = S_1 \oplusdots S_l$ denote the
decomposition of the syzygy bundle into indecomposable locally free sheaves on $C$.
Set $\nu_j = -\deg( S_j) /\rk(S_j) \degc $
Write $\lceil q \nus_j \rceil = q \nus_j + \pers_j (q) $ with
the eventually periodic functions $\pers_j (q)$.
Then the Hilbert-Kunz function of $I$ is
$$ \hkf (q) =e_{HK}(I) q^{2} + \gamma(q) \, ,$$
where
$$e_{HK}(I)  = \frac{1}{2 \degc} \big( \sum_{j=1}^l \frac{\deg(S_j)^2}{\rk (S_j)}
 - \degc^2 \sum_{i=1}^n  d_i^2  \big)$$
is a rational number and where
$$\gamma(q) =  \sum_{j=1}^l \frac{\rk(S_j) \degc}{2} \pers_j (q) (1- \pers_j (q) )
+ \sum_{j=1}^l h^1 (S_j^q( \lceil q \nu_j \rceil )) - n+1 \, .$$
If $q \nu_j$ is not a natural number, then the $h^1$-term is $0$.
In any case the $h^1$-term is eventually constant
and $\gamma(q)$ is an eventually periodic function.
\end{theorem}
\proof
We look at the sequence 
$$ 0 \lra \Syz(f_1 \komdots f_n)(m) \lra
\bigoplus_{i=1}^n  \O(m-d_i) \stacklra {f_1 \komdots f_n}  \O(m) \lra 0 \, ,$$
where $\Syz(f_1 \komdots f_n)(m) \cong S_1(m) \oplusdots S_l (m)$.
The pull-back of this short exact sequence
under the $e$-th absolute Frobenius morphism $F^{e}: C \ra C$ yields
$$ 0 \lra (F^{e}(\Syz(f_1 \komdots f_n))) (m) \lra
\bigoplus_{i=1}^n  \O(m-qd_i) \stacklra {f_1^q \komdots f_n^q}  \O(m) \lra 0 \, ,$$
where 
$F^{e}(\Syz(f_1 \komdots f_n))) (m) \cong S_1^q(m) \oplusdots S_l^q (m)$.
We have  $\hkf (q)= \sum_{m=0}^\infty \length  ((R/I^{[q]}) _m)$ and 
$$  \length  ((R/I^{[q]}) _m)
= h^0 (\O(m)) - \sum_{i=1}^n  h^0(\O(m-qd_i )) + h^0(\Syz (f_1^q \komdots f_n^q)(m) )  \, . $$
We sum these expressions from $m=0$ up to $m =q \nuo$, $\nuo \gg 0$.
We apply Lemma \ref{indecomposablesections} to these sheaves.
Lets denote for fixed numbers
$\nuu \ll 0$ and $ \nuo \gg 0$
the expressions from Lemma \ref{indecomposablesections} for an indecomposable sheaf $T$ by $\expr (T)$. Then
we have to compute $ \expr (\O) - \sum_{i=1}^n \expr (\O(-d_i)) + \sum_{j=1}^l \expr (S_j)$.
We have written down the sum  of the $\deg(T)^2$-terms.
Since the rank and the degree are additive in a short exact sequence, the corresponding terms sum up to $0$.
We have written down the periodicity terms and the $h^1$-term coming from
$T=\Syz(f_1 \komdots f_n)$.
For the invertible sheaves $\O(d)$ there is no periodicity term,
since $\nu =- \deg \O(-d)/ \degc = d$ is always rational,
and the $h^1$-term is always $1$.
This together gives the result.
\qed

\begin{remark}
It is true in general for an $R_+$-primary homogeneous ideal $I$ in a normal two-dimensional
standard-graded domain $R$ that the Hilbert-Kunz function
has the form
$\hkf(q)=e_{HK}(I)q^2 + \gamma(q)$
with a rational Hilbert-Kunz multiplicity $e_{HK}(I)$
and a bounded function $\gamma(q)$, see \cite{brennerhilbertkunzfunction}.
Moreover, if the base field is the algebraic closure of a finite field, then
$\gamma(q)$ is an eventually periodic function. This finiteness condition is not necessary
in the elliptic case; for higher genus this question is open, see also
\cite{hunekemcdermottmonsky}.
\end{remark}

\begin{corollary}
\label{hilbertkunzfunctionsemistable}
Let $k$ denote an algebraically closed field of positive characteristic $p$.
Let $C \subset \PP^N$ denote an elliptic curve
of degree $\degc = \deg (C) = \deg (\O(1)|C)$.
Set
$R= \bigoplus_{n \in \NN} \Gamma(C, \O(n))$.
Let $I$ denote a homogeneous
$R_+$-primary ideal in $R$ given by homogeneous
ideal generators $I=(f_1 \komdots f_n)$ of degree $d_i =\deg(f_i)$
and suppose that the syzygy bundle
$\Syz(f_1 \komdots f_n) $ is semistable.
Then the Hilbert-Kunz function of $I$ is
$$ \hkf (q) =e_{HK}(I) q^{2} + \gamma(q) \, ,$$
where
$$e_{HK}(I)  = \frac{\degc}{2 } \big( \frac{(\sum_{i=1}^n d_i)^2}{n-1} -  \sum_{i=1}^n d_i^2  \big)$$
and
\begin{eqnarray*}
\gamma(q)
&=&   \frac{\degc}{2} \frac{-q \sum_{i=1}^n d_i \modu (n-1) }{n-1}
(1- \frac{-q \sum_{i=1}^n d_i \modu (n-1) }{n-1}  ) \cr
& & + h^1 (\Syz(f_1^q \komdots f_n^q)( \lceil \frac{ q \sum_{i=1}^n d_i}{n-1} \rceil )) - n+1 \, .
\end{eqnarray*}
\end{corollary}
\proof
This follows easily from theorem \ref{hilbertkunzfunctionelliptic} using
the relationships  \ \ 
$$\deg(\Syz(f_1 \komdots f_n)(0))= - \degc \sum_{i=1}^n d_i, \, \,
\nu = \frac{ \degc \sum_{i=1}^n d_i}{(n-1) \degc} =\frac{ \sum_{i=1}^n d_i}{n-1} $$
and
$$\pers (q)=  \frac{( - q \degc \sum_{i=1}^n d_i \modu (n-1) \degc) }{ (n-1) \degc} =
\frac{  -q \sum_{i=1}^n d_i \modu n-1  }{n-1} \, .$$
\qed

\section{The Hilbert-Kunz function of an ideal generated by linear forms}

We bring now the results of the previous sections together
to compute the Hilbert-Kunz function of
an ideal $I \subset R$ generated by an embedding linear system,
where $R= \bigoplus_{n \in \NN} \Gamma(C, \O(n))$
is the graded section ring associated to a very ample invertible sheaf $\O(1)$
on an elliptic curve $C$ over an algebraically closed field $k$ of positive characteristic $p$.
If $C \subset \PP^N=\Proj K[X_0 \komdots X_N]$ is
the corresponding embedding, then
$I=(X_0 \komdots X_N)$ and the syzygy bundle $\Syz(X_0 \komdots X_N) = \Omega_{\PP^N|C}$.
Theorem \ref{hilbertkunzfunctionelliptic} gives at once the following general corollary.
If we want to say more about the behavior of the Hilbert-Kunz function
then we need more information about the restriction of the cotangent bundle.

\begin{corollary}
\label{hilbertkunzfunctionellipticcor}
Let $k$ denote an algebraically closed field of positive characteristic $p$.
Let $C \subset \PP^N= \Proj k[X_0 \komdots X_N]$ denote an elliptic curve
of degree $\degc = \deg (C) = \deg (\O(1)|C)$.
Set $R= \bigoplus_{n \in \NN} \Gamma(C, \O(n))$.
Let $\Syz(X_0 \komdots X_N) = \Omega_{\PP^N}|C= S_1 \oplusdots S_l$ denote the
decomposition of the restricted cotangent bundle into indecomposable locally free sheaves on $C$.
Set $\nu_j = -\deg( S_j) /\rk(S_j) \degc $
Write $\lceil q \nus_j \rceil = q \nus_j + \pers_j (q)$ with
the eventually periodic functions $\pers_j (q)$.
Then the Hilbert-Kunz function of the ideal
$I=(X_0 \komdots X_N) \subset R$ is
$$ \hkf (q) =e_{HK}(I) q^{2} + \gamma(q) \, ,$$
where the Hilbert-Kunz multiplicity is
$$e_{HK}(I)  = \frac{1}{2 \degc} \big( \sum_{j=1}^l \frac{\deg(S_j)^2}{\rk (S_j)}
 - \degc^2 (N+1)  \big)$$
and where
$$\gamma(q) =  \sum_{j=1}^l \frac{\rk(S_j) \degc}{2} \pers_j (q) (1- \pers_j (q) )
+ \sum_{j=1}^l h^1 (S_j^q( \lceil q \nu_j \rceil )) - N \, .$$
is an eventually periodic function.
\end{corollary}
\proof
This follows directly from Theorem \ref{hilbertkunzfunctionelliptic} using $n=N+1$.
\qed

\medskip
If the elliptic curve $C \subset \PP^N$ is embedded by a complete linear system
\-- that is $(X_0 \komdots X_N) =\Gamma(C, \O(1))$ \--
then we know by theorem \ref{tangentrestriction} that the restriction of the cotangent bundle is semistable.
Therefore we obtain the
following result of Fakhruddin and Trivedi (see \cite[Corollary 3.18]{fakhruddintrivedi})
as a corollary.

\begin{theorem}
\label{functioncomplete}
Let $C \subset \PP^N$ denote an elliptic curve over an algebraically closed field $k$
of positive characteristic $p$ embedded by a complete linear system $|\O_C(1)|$.
Then the Hilbert-Kunz function of the maximal ideal in the
homogeneous coordinate ring $R=k[X_0 \komdots X_N]/ \foa$ is given by
$$ \hkf(q) =e_{HK}(R) q^2 + \gamma(q)$$
where $e_{HK}(R)=  \frac{(N+1)^2}{2N} $
and where
$$\gamma(q)= \frac{N+1}{2}(-q \modu N) (1- \frac{(-q \modu N)}{N})
+h^1( \Omega (  \lceil \frac{ q(N+1)}{N}  \rceil) ) -N  \, .$$
If $N \neq p^e$, then the $h^1$-term vanishes.
\end{theorem}
\proof
The curve $C \subset \PP^N$ is projectively normal.
Therefore the domain
$R=\bigoplus_{n \in \NN} \Gamma(C, \O(n))$ is normal and standard-graded, that is generated by finitely many forms of degree one, hence $R \cong k[X_0 \komdots X_N]/ \foa$.
Due to the assumption, the dimension of the global sections of
$ \O_C(1)=\O_{\PP^N}(1)|_C$ is $N+1$.
The global dimension and the degree of
an ample invertible sheaf on an elliptic curve are the same. The degree of $\O_C(1)$ is the degree
of the embedding $C \subset \PP^N$, hence $\deg(C) =N+1$
The restriction of the tangent bundle $\shT_{\PP^N}$ to $C$
is semistable due to theorem \ref{tangentrestriction} and therefore strongly  semistable,
hence we apply corollary \ref{hilbertkunzfunctionellipticcor} with $l=1$.
The rank of the syzygy bundle is $N$ and its degree is $ - (N+1)^2$.
This gives for the Hilbert-Kunz multiplicity $e_{HK}(\fom)=$
$$
= \frac{1}{2 \degc}  \big(\frac{\deg (\Syz) ^2}{\rk (\Syz) } -(N+1)^2 \big)
= \frac{1}{2(N+1)} \big( \frac{ (N+1)^4}{N} - (N+1)^3 \big) 
=\frac{ (N+1)^2}{2N } \, .$$
The $O(q^0)$-term is
$\gamma(q)=  \frac{N(N+1)}{2} \pers (q) (1- \pers (q)) + h^1( \Omega( \lceil q \nu \rceil)) -N $
where $\nu= - \deg (\Omega|C)/\rk (\Omega) \degc = (N+1)/N$.
From $\lceil q \nu \rceil = \lceil q(N+1)/N \rceil = q(N+1)/N + \pers (q)$
it follows that
$\pers (q)= (-q(N+1) \modu N) = (-q \modu N)$.
This gives the result.
\qed

\medskip
In the following corollaries we write down the Hilbert-Kunz functions
in the complete case
for $N=2, 3,4$. For $N=2$, the case of a plane smooth cubic,
we get the following Theorem of Buchweitz-Chen
(see \cite[Theorem 4]{buchweitzchenhilbertkunz}).

\begin{corollary}
Let $C \subset \PP^2$ denote a plane elliptic curve over an algebraically closed field
of characteristic $p \geq 3$. Let $R=k[X,Y,Z]/(F)$ denote its homogeneous coordinate ring.
Then the Hilbert-Kunz function of $R$ is
$$ \hkf (q) = \frac{9}{4}q^2 -\frac{5}{4}\, .$$
\end{corollary}
\proof
Since $q$ is odd we have to insert $(-q  \modu 2 )=1$ in the formula of Theorem
\ref{functioncomplete}. This yields
$$\hkf (q) = \frac{9}{4}q^2 + \frac{3}{2} (1- \frac{1}{2} )-2 
= \frac{9}{4}q^2 -\frac{5}{4} \, .$$
\qed

\medskip
The next corollary treats the case of an elliptic space curve embedded by a complete linear system. The homogeneous coordinate ring is a complete intersection ring $R=K[X_0,X_1,X_2,X_3]/(f,g)$,
where $f$ and $g$ are two quadrics.

\begin{corollary}
Let $C \subset \PP^3$ denote an elliptic curve over an algebraically closed field
of characteristic $p \neq 3$ embedded by a complete linear system in projective space.
Then the Hilbert-Kunz function of the homogeneous coordinate ring $R=k[X_0,X_1,X_2,X_3]/(f,g)$ is
$$ \hkf (q)  = \frac{8}{3} q^2 - \frac{5}{3} \, .$$
\end{corollary}
\proof
The formula in Theorem \ref{functioncomplete} yields
$$ \hkf (q)=\frac{8}{3}q^2 +2\big(( - q \modu 3)(1- \frac{(-q \modu 3)}{3})\big) -3 \, .$$
It does not matter whether $q \modu 3 =1$ or $2$, the big bracket yields always $2/3$.
\qed

\begin{corollary}
Let $C \subset \PP^4$ denote an elliptic curve over an algebraically closed field
of characteristic $p \neq 2$ embedded by a complete linear system.
Then the Hilbert-Kunz function of the homogeneous coordinate ring $R$ is
$$ \hkf(q)  = \frac{25}{8} q^{2} - \frac{17}{8}  \, .$$
\end{corollary}
\proof
The formula in Theorem \ref{functioncomplete} yields
$$ \hkf (q)=\frac{25}{8}q^2 +\frac{5}{2}\big(( -q \modu 4)(1- \frac{(-q \modu 4)}{4})\big) -4 \, .$$
The big bracket yields always $3/4$.
\qed

\medskip
In the previous corollaries we have excluded the case $N =p^{e}$ to avoid the cohomological
term $h^1(\Omega( \lceil \frac{  q(N+1)}{N}\rceil ))$ from theorem \ref{functioncomplete}.
We now deal with the case
$N=p=2$ and we recover the main result of
\cite{monskyellipticchartwo} due to Monsky.
Recall the notion of the Hasse invariant (or $p$-rank)
of an elliptic curve in positive characteristic.
The Hasse-invariant is $0$ or $1$ according whether  the Frobenius morphism
$H^1(C, \O_C) \ra H^1(C, \O_C)$ it the zero map or bijective.

\begin{corollary}
Let $C \subset \PP^2$ denote a plane elliptic curve over an algebraically closed field
of characteristic $2$. Let $R=K[X,Y,Z]/(F)$ denote the homogeneous coordinate ring.
Then the Hilbert-Kunz function of $R$ is
$$\hkf (1) =1,\,\, \hkf (2) =8,\, \, \hkf (2^{e})= \frac{9}{4}2^{2e} \mbox{ for } e \geq 2 $$
for Hasse-invariant $0$ and
$$ \hkf (1) =1,\,\,  \hkf (2) =7, \, \,  \hkf (2^{e})= \frac{9}{4}2^{2e} -1 \mbox{ for } e \geq 2 $$
for Hasse-invariant $1$.
\end{corollary}
\proof
By theorem \ref{functioncomplete}
we have $\hkf(1)=1$ and
$$\hkf(q) = \frac{9}{4} q^2 + h^1(\Syz(X^q,Y^q,Z^q)( \frac{3}{2} q)) -2 \, $$
for $q=2^{e} \geq 2$.
Therefore we have to show for Hasse-invariant $0$ that
$h^1(\Syz(X^q,Y^q,Z^q)( \frac{3}{2} q))= 1$ for $q=2$ and $=2$ for $q \geq 4$,
and for Hasse-invariant $1$ that
$h^1(\Syz(X^q,Y^q,Z^q)( \frac{3}{2} q))= 0$ for $q=2$ and $=1$ for $q \geq 4$.
The bundle $\Syz(X^q,Y^q,Z^q)( \frac{3}{2} q)$ has rank two and
its determinant is $\O_C$, therefore it is isomorphic to its dual,
and we can replace $h^1$ by $h^0$.

Up to an automorphism we may assume that $C$ is given by an equation
$X^3+Y^3+Z^3+ \lambda XYZ=0$. Due to \cite[Proposition IV.4.21]{haralg}
the curve has Hasse invariant $0$ if and only if $\lambda =0$.
Moreover, for $\lambda^3=1$ it is easy to see that the point
$(1, \lambda, \lambda)$ is a singular point on the curve, hence $\lambda^3 \neq 1$.

First assume that $\lambda =0$, so we are dealing with the Fermat cubic $X^3+Y^3+Z^3=0$.
This curve equation yields at once the syzygy $(X,Y,Z)$ for $(X^2,Y^2,Z^2)$
of degree $3$
and hence the short exact sequence
$$0 \lra \O_C \lra \Syz(X^2,Y^2,Z^2)(3) \lra \O_C \lra 0 \, .$$
This shows that the first
Frobenius pull-back of the restricted tangent bundle is not stable (though semistable).
It is clear that $h^0(C, \Syz(X^2,Y^2,Z^2)(3))=1$. Therefore this sequence does not split
and it is given by a non-trivial cohomology class $c \in H^1(C, \O_C)$.
Since the Hasse-invariant is $0$,
the Frobenius pull-back of the short exact sequence splits and hence
$$\Syz (X^q,Y^q,Z^q) (m) \cong \O(m-3q/2) \oplus \O(m- 3q/2)$$
for $q \geq 4$.
Therefore the correction term is $h^0(C, \Syz(X^q,Y^q,Z^q)(3q/2))=2$.

We consider now the case of Hasse-invariant is $1$, which means that $\lambda \neq 0$,
$\lambda^3 \neq 1$. It is clear that $H^0(C,\Syz(X^2,Y^2,Z^2)(3)) =0$.
So we look at the syzygies $H^0(C, \Syz(X^4,Y^4,Z^4)(6))$.
Multiplying the curve equation with $X^3+Y^3+Z^3$ yields
\begin{eqnarray*}
& & (X^3+Y^3+Z^3)(X^3+Y^3+Z^3+\lambda XYZ) \cr
&=& X^6+Y^6+Z^6 +\lambda X^4YZ + \lambda XY^4Z + \lambda XYZ^4 \cr
&=& X^4(X^2+ \lambda YZ) + Y^4(Y^2+ \lambda XZ) +Z^4(Z^2+ \lambda XY) 
\end{eqnarray*}
hence $(X^2+ \lambda YZ, Y^2+ \lambda XZ, Z^2+ \lambda XY)$ is a global syzygy
for $(X^4,Y^4,Z^4)$ of degree $6$.
The corresponding short exact sequence is
$$0 \lra \O_C \lra \Syz(X^4,Y^4,Z^4)(6) \lra \O_C \lra 0 \, .$$

We claim that this sequence does not split. To see this we show that
the space of global syzygyies of degree $6$ for $X^4,Y^4,Z^4$ is one-dimensional.
Such a syzygy is given by a multiple of the curve equation fulfilling
$$H(X^3+Y^3+Z^3+ \lambda XYZ) \in (X^4,Y^4,Z^4) \, ,$$
where $\deg H =3$.
The polynomial $H$ is a linear combination of the $10$ monomials of degree $3$.
The monomial $XYZ$ cannot occur, since the product $\lambda X^2Y^2Z^2$
cannot arise in another way.
The products with $X^3,Y^3,Z^3$ yield $X^3Y^3$ etc., and these monomials
do not occur in another way.
If $H=aX^3+bY^3+cZ^3$ we arrive at the syzygy from above.

Hence we may assume that $H$ is a linear combination of monomials of type $(1,2)$,
say $H=aX^2Y+bXY^2+cXZ^2+dX^2Z+eYZ^2+fY^2Z$.
Multiplying with the curve equation yields
(modulo $(X^4,Y^4,Z^4)$)
the terms
$$ (b+e \lambda)XY^2Z^3  +(c+f \lambda)XY^3Z^2  + (d+b \lambda)X^2Y^3Z $$
$$+(a+c \lambda)X^2YZ^3    + (e+d \lambda)X^3YZ^2 + (f+a \lambda)X^3Y^2Z    $$
The conditions $b+e \lambda =0$, $d+b\lambda =0$ and $e+d \lambda=0$  yield $e=e \lambda^3$
(and $f =f \lambda^3$).
Since $\lambda^3 \neq 1$ we see $e=f=0$ and hence all coefficients are $0$.
This proves the claim.

Since the Hasse-invariant is $\neq 0$ it follows that also all the Frobenius pull-backs
of the short exact sequences do not split.
From this it follows that $h^0(C, \Syz(X^q,Y^q,Z^q)(3q/2))=1$ for $q \geq 4$.
\qed

\section{Elliptic space curves}

In this last section we deal with the Hilbert-Kunz function of an elliptic space curve
$C \subset \PP^3$. The degree of $C$ may be every number $\degc \geq 4$.
For $\degc =4$ the curve is embedded by a complete linear system and the
associated graded section ring
$R= \bigoplus_{n \in \NN} \Gamma(C, \O_C(1))$
is normal and standard-graded.
This is not true for $\degc \geq 5$
and the following theorem computes the Hilbert-Kunz function
of the ideal $(X_0,X_1,X_2,X_3)$ of the graded section ring $R$,
not of the homogeneous coordinate ring.
We call this the Hilbert-Kunz function of the curve $C \subset \PP^3$
embedded by the linear system
$(X_0,X_1,X_2,X_3)= \Gamma(\PP^3, \O(1))|_C \subseteq \Gamma(C, \O(1))$.
The classification from section 3 yields the following list of Hilbert-Kunz functions.

\begin{theorem}
Let $C \subset \PP^3$ denote an elliptic space curve of degree $\delta$.
Then there exist the following three possibilities for the Hilbert-Kunz function of $C$,
depending on the splitting behavior of $\Omega_{\PP^3}|_C$.

\numiii

\begin{enumerate}

\item
Suppose that $\Omega|C$ is indecomposable.
Then
$$\hkf (q)= \frac{2 \degc}{3} q^2  + \frac{\degc}{3}  -3  \, .$$
for $p \neq 3$
and
$$ \hkf (q)= \frac{2 \degc}{3} q^2 
+ h^1 ( \Syz(X_0^q,X_1^q,X_2^q,X_3^q)(  \frac{4q}{3} )) -3 $$
for $p=3$ (for $q \neq 1$).

\item
Suppose that
$\Omega|C \cong S \oplus M$, where $S$ is indecomposable of rank two and $M$ is invertible.
Set $\nu_1 =- \deg (S)/2\degc$ and $\nu_2 = - \deg(M)/\degc$
and write $\lceil q \nu_j \rceil = q \nu_j + \pers_j (q)$, $j=1,2$.
Then
\begin{eqnarray*}
\hkf (q)\!\!\!\! &=&\!\!\!\! \frac{1}{2 \degc} ( \frac{ \deg (S)^2}{2}\! +\! \deg(M)^2 \!-\!4 \degc^2 )   q^2 
\!+ \!\degc  \pers_1 (q)\! (1- \pers_1 (q))\! +\! \frac{\degc}{2} \pers_2 (q)\!(1 - \pers_2(q))
\cr
& & + h^1 (S^q( \lceil \frac{ -q \deg(S)}{2 \degc} \rceil)) 
+ h^1 (M^q( \lceil \frac{ -q \deg(M)}{\degc} \rceil)) 
-3  
\end{eqnarray*}

\item
Suppose that
$\Omega|C \cong M_1 \oplus  M_2 \oplus M_3$ with three invertible sheaves $M_1$, $M_2$ and $M_3$.
Set $\nu_j = - \deg (M_j)/\degc$ and write $\lceil q \nu_j \rceil = q \nu_j + \pers_j(q)$, $j=1,2,3$.
Then
\begin{eqnarray*}
\hkf (q)\!\!\! &=& \!\!\! \frac{ \deg (M_1)^2 + \deg(M_2)^2 + \deg(M_3)^2 -4 \degc^2}{2 \degc}   q^2 
+ \frac{\degc}{2}( \sum_{j=1}^3  \pers_j(q) (1 - \pers_j(q))) \cr
& & + h^1 (M_1^q( \lceil \frac{ -q \deg(M_1)}{\degc} \rceil)) 
+ h^1(M_2^q( \lceil \frac{ -q \deg(M_2)}{\degc} \rceil)) \cr
& & + h^1 (M_3^q( \lceil \frac{ -q \deg(M_3)}{\degc} \rceil)) -3  
\end{eqnarray*}

\end{enumerate}

\end{theorem}

\proof
The rank of $\Omega|C$ is $3$ and its degree is $-4 \degc$.
We compute the Hilbert-Kunz multiplicity as 
$e_{HK} = \frac{1}{2 \degc} ( \sum_{j} \frac {\deg (S_j)^2}{\rk (S_j)} - 4 \degc^2 )$
(theorem \ref{hilbertkunzfunctionelliptic}).
This gives in the indecomposable case
$\frac{1}{2 \degc} (  \frac {16 \degc ^2}{3} - 4 \degc^2 )= \degc (\frac{8}{3} -2) =2 \degc/3 $.
In the other cases we just have inserted.

Now we look at the lower term $\gamma(q)$ coming from theorem \ref{hilbertkunzfunctionelliptic}.
In the indecomposable case
we have by remark \ref{periodicityremark}
$$\pers (q) = (q \deg (S)  \modu 3 \degc)/ 3 \degc = ( -4q \modu 3)/3 =(-q \modu 3)/3 \, .$$
For $p \neq 3$ we have $\frac{(-q \modu 3)}{3}  (1- \frac{(-q \modu 3)}{3}) = \frac{2}{9}$,
for $p=3$ this term is $0$ for $q \neq 1$.
In cases (ii) and (iii) we have just written down the formula for $\gamma(q)$
coming from theorem \ref{hilbertkunzfunctionelliptic}.
\qed

\begin{remark}
Suppose that $C \subset \PP^3$ is an elliptic space curve of degree $\degc$.
If $\Omega|C$ is semistable, then the Hilbert-Kunz multiplicity is $2 \degc/3$.
This is automatically true if $\Omega|C$ is indecomposable,
but this holds also in case (ii) if
$\deg (M)= \deg(S)/2$ and in case (iii) if the degrees of
the invertible sheaves $M_1$, $M_2$ and $M_3$ are equal.
If $\Omega|C$ is not semistable, then $e_{HK}(C) > 2 \degc/3$.
This follows from a direct computation or from
\cite[Lemma 1.3]{brennerhilbertkunzcriterion}.
\end{remark}

\bibliographystyle{plain}

\bibliography{bibliothek}

\end{document}